\documentclass[11pt]{article}
\usepackage{amsfonts}

\usepackage{graphics}
\usepackage{indentfirst}
\usepackage{cite}
\usepackage{latexsym}
\usepackage{amsmath}
\usepackage{amssymb}
\usepackage[dvips]{epsfig}
\usepackage{amscd}

\newtheorem{theorem}{Theorem}[section]
\newtheorem{remark}{Remark}[section]
\newtheorem{lem}{Lemma}[section]

\newtheorem{lemma}[theorem]{Lemma}

\newcommand{\n}{\rho}

\def\pf{{\it Proof.}  }

\newcommand{\pa}{\partial}

\newcommand{\bi}{\bibitem}

\newcommand{\bt}{\begin{theorem}}
\newcommand{\bl}{\begin{lemma}}
\newcommand{\el}{\end{lemma}}
\newcommand{\et}{\end{theorem}}
\newcommand{\ga}{\gamma}

\newcommand{\al}{\alpha }
\newcommand{\de}{\delta}
\newcommand{\ve}{\varepsilon}
\newcommand{\la}{\label}
\newcommand{\si}{\sigma}

\newcommand{\bn}{\begin{eqnarray}}
\newcommand{\en}{\end{eqnarray}}
\newcommand{\bnn}{\begin{eqnarray*}}
\newcommand{\enn}{\end{eqnarray*}}

\newcommand{\bnnn}{\begin{eqnarray*}}
\newcommand{\ennn}{\end{eqnarray*}}
\newcommand{\ben}{\begin{enumerate}}
\newcommand{\een}{\end{enumerate}}

\newcommand{\ba}{\begin{aligned}}
\newcommand{\ea}{\end{aligned}}
\newcommand{\be}{\begin{equation}}
\newcommand{\ee}{\end{equation}}

\def\norm[#1]#2{\|#2\|_{#1}}

\def\xix{\int_0^1}

\def\xiT{\int_0^T}

\def\O{\Omega}

\def\xl{\left}
\def\xr{\right}

\makeatletter      
\@addtoreset{equation}{section}
\makeatother       

\title{On Global
Classical Solutions to 1D Compressible  Navier-Stokes Equations with Density-Dependent Viscosity and Vacuum\thanks{B. L\"u is supported by NNSFC (Nos. 11601218 \& 11771382), Science and Technology Project  of Jiangxi Provincial Education Department (No. GJJ160719).}}
\date{}
\author{ Boqiang L\"u\thanks{College of Mathematics and Information Science, Nanchang Hangkong University, Nanchang 330063, P. R. China({\tt lvbq86@163.com}). }
 \quad Yixuan Wang\thanks{Institute of Applied Mathematics, AMSS,
Chinese Academy of Sciences, Beijing 100190,  P. R. China ({\tt wangyixuan\_14@163.com}).
 } \quad Yuhang Wu\thanks{Department of Physics, Beijing Normal University, Beijing, 100875, P. R. China
 ({\tt 1095986995@qq.com}).}}

\begin{document}
\maketitle

\begin{abstract} For the  initial boundary value problem of   compressible barotropic Navier-Stokes equations in one-dimensional bounded domains with   general density-dependent viscosity and large  external force,  we prove that there exists   a unique  global  classical  solution with large initial data containing vacuum.  Furthermore, we show that the density is bounded from above independently of time which in particular yields the large time behavior of the solution as time tends to infinity: the density and the velocity converge  to the steady states in $L^p$ and  in $W^{1,p}$  ($1\le p<+\infty$) respectively. Moreover, the decay rate in time of the solution is shown to be exponential.
Finally, we also prove  that the spatial gradient of the density
 will blow up as time tends to infinity when vacuum states appear initially even at one point.
\end{abstract}


\section{Introduction and main results}
  We consider the one-dimensional
compressible  Navier-Stokes equations which read as follows:
\be\la{1d}
\begin{cases} \rho_t + (\rho u)_x = 0,\\
 (\rho u)_t + (\rho u^2)_x + [P(\n)]_x -[\mu(\n)u_x]_x=\n f.
\end{cases}
\ee
Here $t\ge 0$ is time, $x\in \O=(0,1)$ is the spatial coordinate,
 $\n$ and $ u $  represent respectively the fluid density and velocity.
The pressure $P$ is given by
\be\la{n2}
P(\rho) = A\rho^{\gamma}\, ( A>0,\, \ga>1).
\ee
 In the sequel, without loss of generality, we set $A= 1.$ The viscosity   $\mu(\n)$   satisfies
\be\la{1d1}
~~0<\bar{\mu}\le \mu(\n),~~\forall~\n\ge 0,
\ee where $\bar{\mu}$ is  a given positive constant. The function $f=f(x)$ is the external force.
 We look for the solutions, $(\n(x,t),u(x,t)),$ to
 the initial-boundary-value problem with  the    boundary conditions:
\be\la{1d2}
u(0,t)=u(1,t)=0,~~~t\ge 0,
\ee  and the initial ones:
\be \la{n3}
(\n,\n u)(x,0)=(\n_0 ,\n_0u_0)(x ).
\ee

There is huge literature on the studies of the global existence and large time behavior of solutions to the compressible Navier-Stokes equations.
For the initial density away from vacuum, there are  many results concerning the global existence and large-time dynamics of solutions to the
 one-dimensional(1D) problem, see \cite{Kaz,serre1,serre2,Hof,v1989,is2002,is2003,lij1d}  and the reference therein.
%
%
For the case of density-dependent viscosity,       Liu-Xin-Yang \cite{lxy} obtains that  the viscosity of a gas depends on the temperature for the non-isentropic case, and thus on
the density for the isentropic case. For $\mu(\n)=\n^\theta$, under different restrictions
on the index $\theta$ and the regularities  of initial data,  the global existence  of solutions to 1D compressible Navier-Stokes
equations is investigated  in \cite{ok,jxz,yyz,qyz,fz,yz,ljcmp2008} and the references therein.
%
 When $\mu(\n)$  admits a positive constant lower bound,
 the global well-posedness of solutions without initial vacuum to 1D problem is discussed extensively (see \cite{kawohl, zl89,v1989,is2002,is2003}  and  the references therein).
 Recently, in addition to \eqref{1d1},  under some additional stringent condition on $\mu\in C^2[0,\infty)$:
\be\la{wen}
  \mu(\n)\le C(1+P(\n)) ~~\forall~\n\ge 0,
\ee Ding-Wen-Zhu \cite{wen2011} proves the  global existence of classical large solutions to   \eqref{1d}-\eqref{n2} with vacuum. However, since the   upper bound of the density obtained in   \cite{wen2011}  depends crucially on   time,   nothing is known concerning the large-time behavior of the solutions in \cite{wen2011}. For the initial density away from vacuum,  Stra\~skraba-Zlotnik \cite{is2002,is2003}     proves the large time behavior of the solution $(\n,u)$ to \eqref{1d}.  More precisely, they show that the solution $(\n,u)$ tends to the stationary case $(\n_s,0)$ as time tends to infinity. Here,  the  stationary density $\n_s$  is a solution to the  stationary problem:
  \be\la{1s}
\begin{cases}  [P(\n_s)]_x=\n_s f,~~~x\in (0,1),\\
\xix\n_sdx=\xix\n_0dx.\end{cases}
\ee

In this paper, we will study  the global existence and large-time behavior of  strong and  classical solutions to \eqref{1d}-\eqref{n3} not only for   general density-dependent  viscosity \eqref{n3} which is independent of $P(\n)$ but also for the density containing vacuum initially.   

Before stating the main results, we first explain the notations and
conventions used throughout this paper.
For   $1\le r\le \infty, k\ge 1$,
$$L^r=L^r(0,1 ),\quad W^{k,r}  = W^{k,r}(0,1) , \quad H^k = W^{k,2}(0,1).$$
Moreover,  without loss of generality, assume that the initial density $\n_0$ satisfies
 \be\la{oy3.7} \xix  \n_0dx=1.\ee

Our first result concerns the global existence of strong solutions with large initial data.

 \begin{theorem}\la{th11}   Suppose that $f\in H^1$, $\mu \in C^1[0,\infty),$ and the initial data $(\n_0,u_0)$ satisfy
\be\la{1d3a}
  0\le\n_0\in H^1,  ~u_0\in   H^1_0. \ee  Then, there exists a unique strong  solution $(\n,u)$ to the  problem  \eqref{1d}-\eqref{n3} satisfying for any $0< T<\infty,$
   \be
   \la{1d5}\begin{cases}
   \rho  \in C([0,T]; H^1)\cap H^1(0,T;L^2), \\
    u\in L^\infty(0,T; H^1_0) \cap L^2(0,T;H^2), \\t^{1/2}u \in L^\infty(0,T; H^2),
~~ t^{1/2}u_t\in L^2(0,T; H^1_0)  .\end{cases} \ee
 Moreover,  the density remains uniformly  bounded for all time, that is,
 \be\la{1d6}
  \sup_{0\le t<\infty}\|\n(\cdot,t)\|_{L^\infty} <\infty,
\ee  and the following large-time behavior holds:
\be \la{1.12'}\lim_{t\rightarrow \infty} \| u(\cdot,t)\|_{W^{1,p}} =0,~~~~\forall~p\in [1,\infty).\ee
\end{theorem}

Then the following result shows that the strong solutions obtain by Theorem \ref{th11} become classical provided initial data $(\n_0,u_0)$ satisfy some additional conditions.
  \begin{theorem}\la{th1}  In addition to \eqref{1d3a}, suppose that $f\in H^1$, $\mu \in C^2[0,\infty),$  and the initial data $(\n_0,u_0)$ satisfy
\be\la{1d3}
 \n_0\in H^2, ~P(\n_0)\in H^2,~u_0\in H^2\cap H^1_0, \ee
and the following compatibility condition:
 \be\la{1d4}
[\mu(\n_0)u_{0x}]_x-[P(\n_0)]_x=\n_0^{1/2}g ,  \ee
for some $g\in L^2$. Then, the strong solution $(\n,u)$ obtained by Theorem \ref{th11} becomes  classical and satisfies for any $0< T<\infty,$
   \be
   \la{1d5'}\begin{cases}
   \rho,   P(\n) \in C([0,T]; H^2), \\
    u\in C([0,T]; H^2) \cap L^2(0,T;H^3), \\
 u_t\in L^2(0,T; H_0^1),~~t^{1/2}u \in L^\infty(0,T; H^3),\\
t^{1/2}u_t \in L^\infty(0,T; H^1)\cap L^2(0,T; H^2),~~  t^{1/2}\n^{1/2}u_{tt}\in L^2(0,T; L^2).\end{cases} \ee

\end{theorem}

For further  studying the large-time behavior of the strong solutions,   we  first state some known results about the existence and uniqueness of  positive solutions to the stationary problem \eqref{1s}, which has been  discussed extensively under different conditions (see for example \cite{zl1,zaa95,is2002,is2003}).
 \begin{lem}[\cite{zaa95}]\la{le2}   If $f\in L^\infty$ satisfies
 \be\la{niu}\xix \xl(1- \gamma^{-1}\xr)^{\frac{1}{\gamma-1}} \xl(\int_0^x f(y)dy-\min_{[0,1]}\int_0^x f(y)dy\xr)^{\frac{1}{\gamma-1}} dx<\int_0^1\rho_0dx,\ee there exists a unique positive solution $\n_s$ to \eqref{1s} which satisfies
 \be\la{pre1}
\n_s\in W^{1,\infty},~~  0<K_1\le \n_s\le K_2,
\ee
where $K_1$ and $K_2$  are positive constants depending on $\|f\|_{L^\infty}$.
\end{lem}

Then, we have the following result concerning  the  large time asymptotic behavior  of  strong solutions.

\begin{theorem}\la{th2}   Under the same conditions as in Theorem \ref{th11}, for  $f$ satisfying \eqref{niu},   there are    positive constants $\al$ and $C$ depending only on the initial data and $\|f\|_{H^1}$  such that for any $p\in [1,\infty)$ and any $t>1,$
\be\label{ed1}
 \|\n(\cdot,t)-\n_s(\cdot)\|_{L^p}+ \|u(\cdot,t)\|_{W^{1,p}}\le C  e^{-\frac{\al}{ p} t}.
\ee
 Moreover, if  there exists some point $x_0\in [0,1]$ such that $\n_0(x_0)=0,$  then the spatial gradient of the density the unique strong solution $(\n, u)$ to  the  problem  \eqref{1d}-\eqref{n3}  has to blow up as $t\rightarrow \infty$ in the following sense,
\be\label{bp}
\lim_{t\rightarrow \infty}\|\n_x(\cdot,t)\|_{L^r}=\infty,~~~\forall~r\in (1,\infty).
\ee
\end{theorem}



A few remarks are in order:

\begin{remark}\la{re1}
It should be noted here that the solution  $(\n,u)$  obtained in Theorem \ref{th1} is actually a classical one to  \eqref{1d}-\eqref{n3}.
Indeed, by the Sobolev embedding theorems, we have
$$H^k(0,1)\hookrightarrow C^{k-\frac{1}{2}}[0,1],~~~\mbox{for}~k=1,\cdots,3, $$
which together with   \eqref{1d5'}  gives
\be \label{ccc}(\n,\,P,\, u) \in C([0 ,T];C^{1+\frac{1}{2}}[0,1]),~~ \n_t   \in C([0 ,T];C^{\frac{1}{2}}[0,1]).\ee
Furthermore, one can deduce from \eqref{1d5'} that  for any $0<\tau<T$,
$$u \in L^\infty(\tau,T; H^3),~~~u_t \in L^\infty(\tau,T; H^1)\cap L^2(\tau,T; H^2),$$
which yields that for $0<\al<1/2$,
\be \label{ccc1}u \in   C([\tau,T]; C^{2,\alpha}[0,1]),~~~u_t \in   C([\tau,T]; C^{\alpha}[0,1]).\ee
 Hence, it follows from \eqref{ccc} and \eqref{ccc1}   that  $(\n,u)$   is  a classical solution to  \eqref{1d}-\eqref{n3}.
\end{remark}

\begin{remark} \la{re4} To obtain the global  existence  of strong solutions in Theorem \ref{th11}, we do not need the additional  compatibility condition \eqref{1d4}. Indeed, we only need the initial data satisfying  the compatibility condition \eqref{1d4} for some $g\in L^2$     in proving the global well-posedness of classical solution in our Theorem \ref{th1}, which is in sharp contrast to \cite{wen2011} where they need   $g\in H^1.$ Therefore, our theorems  essentially weaken those assumptions on the compatibility condition in \cite{wen2011}.
\end{remark}

\begin{remark} \la{re2}
In our result, the only restriction  on $\mu\in C^2[0,\infty)$ is  \eqref{1d1}, which is much more general than those in \cite{wen2011} where   in addition to \eqref{1d1},   they need an additional stringent condition   \eqref{wen} which  plays a crucial role in the analysis of \cite{wen2011}.
\end{remark}

\begin{remark}\la{re3} In Theorem \ref{th11},
we obtain the time-independent upper bound of the density  in \eqref{1d6} and  the large-time behavior of  the velocity \eqref{1.12'}, which are in sharp contrast to \cite{wen2011} where  the corresponding a priori ones  depend on time.
Moreover, these results also generalized the similar ones in \cite{is2002,is2003} where they need   initial density strictly away from vacuum to the case that the density allows vacuum intially.
\end{remark}

\begin{remark}\la{re5}
In contrast to  the results in \cite{is2002,is2003}, where     the spatial $L^2$-norm  $\n_x$ is  proved to be bounded independently of time provided the density strictly away from vacuum, our Theorem \ref{th2}  shows that  the spatial $L^r$-norm ($r>1$) of $\n_x$  will blow up when time goes to infinity provided  vacuum appears initially. 
\end{remark}

We now make some comments on the analysis of this paper. We begin with the local existence theorem (see  \cite{coi1} or Lemma \ref{le3} below) of classical solutions
to problem \eqref{1d}-\eqref{n3} with  the initial density  strictly away from vacuum.
Then, we prove that the local strong (classical) solution with vacuum also exists and thus extend the local existence time to be global. Hence, we need some
global a priori estimates which are independent of the lower bound of density. It turns out that
the key issue is to derive both the time-independent lower order estimates  and
the time-dependent higher order ones (see Section 2).
It should be noted  that the methods used in  Ding-Wen-Zhu \cite{wen2011} can not be adapted here.
Indeed, on the one hand, the analysis in \cite{wen2011} relies heavily on the special assumption on  viscosity $\mu(\n)$ (see \eqref{wen}), that is $\mu(\n)$ should be bounded by   $P(\n)$ pointwisely.
On the other hand, it seems difficult to study  the large-time behavior of solutions since the a priori  estimates obtained in \cite{wen2011} are all   time-depending. To overcome these difficulties,
motivated by   Li-Xin \cite{lijde,lx1},
we succeed in obtaining  the  key uniform upper bound of the density by making full use
of Zlotnik  inequality (see Lemma \ref{la3.2}),  and bounding the $L^2$-norm of $u_x$ according to the  material derivative  $\dot u$ (see Lemma \ref{la3.5}).
The time-dependent higher order estimates of $(\n, u)$ are derived by some standard arguments and
the time-weighted estimates due to \cite{H3} (see also \cite{hlx1,lx1,lisiam}).
Next, with both the uniform upper bound of the density  and the time-independent lower order estimates at hand,
 we use the methods owing to Stra\v{s}kraba-Zlotnik \cite{is2002,is2003} and
 thus prove  the following large-time behavior
  $$  \lim_{t\rightarrow \infty}\xl(\| u(\cdot,t)\|_{W^{1,p}}+\| \n-\n_s\|_{L^{p}}\xr)=0,~~~~\forall~p\in [1,\infty). $$
Finally, using a key testing function $\n_s^{-1}\int_0^x(\n-\n_s)dy$
 motivated by Huang-Li-Xin \cite{jmpa06} and Li-Zhang-Zhao \cite{lisiam},
  we derive the desired  exponential decay rate estimate \eqref{ed1}
  in Theorem \ref{th2}(see Section 4).

The rest of the paper is organized as follows: In Section 2, we will  derive the necessary
a priori estimates on smooth solutions. The main results,
Theorems \ref{th11}--\ref{th1}    and \ref{th2},  are proved in Sections 3   and   4 respectively. 

 \section{A priori estimates}

In this section, we will establish some necessary a priori bounds
for smooth solutions to the problem \eqref{1d}-\eqref{n3} to extend the local classical solution guaranteed by following
Lemma \ref{le3}, whose proof can be completed by similar arguments as in  \cite{coi1,HLM}.
 \begin{lemma}\la{le3}   Assume that $f\in H^2$ and  the initial data $(\n_0, u_0)$ satisfies
   \bnn\la{pre4}
0<\de\le\n_0,~~ \n_0\in H^3, ~P(\n_0)\in H^3,~u_0\in H^3\cap H^1_0.
\enn
Then, there exists a small time $T_0>0$ and a unique classical solution $(\n,u)$ to the problem \eqref{1d}-\eqref{n3}  on $(0,1)\times(0, T_0]$  such that
   \bnnn
   \la{pre3}\begin{cases}
   \rho,   P(\n) \in C([0,T_0]; H^3), \\
 u\in C([0,T_0]; H^3\cap H_0^1) \cap L^2(0,T_0;H^4),\\
 \rho_t \in C([0,T_0]; H^2),~~
   u_t\in  C([0,T_0]; H^1_0)\cap L^2(0,T_0;H^2).\end{cases} \enn
\end{lemma}

Let $T>0$ be a fixed time and $(\n,u)$ be
the smooth solution to  problem \eqref{1d}-\eqref{n3} with the initial data  $(\n_0,u_0)$ satisfying the  assumptions in Lemma \ref{le3}.

 \subsection{\la{se3} A priori estimates(I): lower order estimates}

This subsection is concerned with the time-independent lower order estimates of  the solution $(\n,u)$.  In this subsection, we will use the convention that $C$ denotes a generic positive constant depending  on
 $\gamma, \bar\mu, \|\n_0\|_{H^1}, \|u_0\|_{H^1}$, and $ \|f\|_{H^1}$
but independent of $T$, and  use $C(a)$ to emphasize that $C$ depends on $a.$

First, we state the following Zlotnik  inequality, whose proof can be found in \cite{zl1}, will be used to get the
uniform (in time) upper bound of the density.
\begin{lemma}[\cite{zl1}]\la{le1}   Let the function $y$ satisfy
\bnn y'(t)\le  g(y)+b'(t) \mbox{  on  } [0,T] ,\quad y(0)=y^0, \enn
with $ g\in C(R)$ and $y,b\in W^{1,1}(0,T).$ If $g(\infty)=-\infty$
and \bnn\la{a100} b(t_2) -b(t_1) \le N_0 +N_1(t_2-t_1)\enn for all
$0\le t_1<t_2\le T$
  with some $N_0\ge 0$ and $N_1\ge 0,$ then
\bnn y(t)\le \max\left\{y^0,\overline{\zeta} \right\}+N_0<\infty
\mbox{ on
 } [0,T],
\enn
  where $\overline{\zeta} $ is a constant such
that \bnn\la{a101} g(\zeta)\le -N_1 \quad\mbox{ for }\quad \zeta\ge \overline{\zeta}.\enn
\end{lemma}

Then,  we will give the key time-independent upper bound for the density as follows.
 \begin{lemma} \la{la3.2}
There is a positive constant $C$   such that for any $(x,t)\in [0,1]\times[0,T],$
  \be\la{gj3.2}
0\le \n(x,t)\le C .
\ee
\end{lemma}
\pf  First, standard energy estimate leads to
 \be\la{gj3.1} \ba
&\sup_{0\le t\le T} \xix \left(\frac{1}{2}\n u^2+\frac{1}{\gamma-1}P\right)dx+\int_0^T\xix\mu(\n)u_x^2dxds\\
&\le C+\sup_{0\le t\le T} \xix \n \int_0^xfdydx\\
&\le C+\|f\|_{L^\infty}\sup_{0\le  t\le T} \xix \n  dx\\
&\le C,
\ea\ee
where in the last inequality one has used the following fact:
 \be\ba\la{03.2}
& \xix \n  dx=\xix\n_0dx=1
\ea\ee
owing to  \eqref{oy3.7} and $ \eqref{1d}_1$.

Next, integrating $\eqref{1d}_2$ over $(0,x)$ gives
  \be\la{l31}\ba
 -\mu(\n)u_x=&-P(\n)+\int_0^x\n fdy -\frac{\pa}{\pa t}\int_0^x\n udy-\n u^2\\&+\left(P(\n)+\n u^2-\mu(\n)u_x\right)(0,t),
\ea\ee
which in particular implies  \be\la{l32}\ba
&\left(P(\n)+\n u^2-\mu(\n)u_x\right)(0,t)\\&=-\xix\mu(\n)u_xdx+\xix P(\n)dx-\xix \int_0^x\n fdydx \\
&\quad+ \frac{d}{dt}\xix\int_0^x\n udydx+\xix \n u^2dx.
\ea\ee
Combining this with   \eqref{l31} shows that for $D_t\triangleq \frac{\pa}{\pa t} +u  \frac{\pa}{\pa x}$
  \be\la{l3}\ba
&-\mu(\n)u_x+P(\n) \\
&= \int_0^x\n fdy-\xix\int_0^x\n fdydx+\int_0^1\n u^2dx+\int_0^1P(\n)dx\\
&\quad+D_t\left(\int_0^1\n\int_2^\n \mu(s)s^{-2}dsdx+\int_0^1\int_0^x\n udydx-\int_0^x\n udy\right)\\
&\le  C+D_t\left(b_1(t)+b_2(t)+b_3(t)\right),
\ea\ee
where one has used \eqref{gj3.1}  and the following fact
  \bnn\ba
-\int_0^1\mu(\n) u_xdx=D_t\int_0^1\n\int_2^\n \mu(s)s^{-2}dsdx
\ea\enn
due to $\eqref{1d}_1$.

Next,  it follows from $\eqref{1d}_1$ that
  \bnn-\mu(\n)u_x=
D_t\int_1^\n \mu(s)s^{-1}ds ,
\enn which together with
  \eqref{l3}   gives
  \be\la{l5}\ba
D_t\int_1^\n \mu(s)s^{-1}ds\le -P(\n)+C+D_t\left(b_1(t)+b_2(t)+b_3(t)\right).
\ea\ee

 Then, on the one hand, we have
\be\la{l8} \ba
b_1(t) \le 2^{-1}\int_1^{\max\{\sup\limits_{\Omega_T}\n,2\}}\mu(s)s^{-1}ds,
\ea\ee
where $\Omega_T=[0,1]\times[0, T]$.
On the other hand, one deduces from \eqref{gj3.1} that
\be\la{l9} \ba
|b_2(t)|+|b_3(t)|&\le C\int_0^1\int_0^x\n dy\int_0^x\n u^2dydx+\int_0^1\n u^2dx\int_0^1\n dx \le C.
\ea\ee
Finally, applying the Zlotnik  inequality (Lemma \ref{le1}) to \eqref{l5}, we get after using  \eqref{l8} and \eqref{l9} that
  \bnn\ba
\int_1^\n \mu(s)s^{-1}ds\le C+2^{-1}\int_1^{\max\{\sup\limits_{\Omega_T}\n,2\}}\mu(s)s^{-1}ds,
\ea\enn
which together with \eqref{1d1}     implies \eqref{gj3.2} and completes the proof of Lemma \ref{la3.2}. \hfill $\Box$

\begin{lemma} \la{la3.5}
Let $(\n,u)$ be a smooth solution of
 \eqref{1d}-\eqref{n3}  on $[0,1] \times [0,T].$ Then
 \be
\ba\la{gj3.5}
   \sup_{0\le t\le T}  \| u \|_{H^1}^2+\xiT \|\n^{1/2}\dot u\|_{L^2}^2 dt\le C,
\ea
\ee where and in what follows, we denote
$ \dot u\triangleq u_t+uu_x. $
\end{lemma}

\pf First, multiplying $\eqref{1d1}_{2}$  by  $ \dot u$ and   integrating the resulting equation by parts yield
\be
\ba\la{3.2}
  &\frac{1}{2}\frac{d}{dt}\xix \mu(\n)u_x^2dx+\xix\n\dot u^2dx\\
  &=\frac{d}{dt}\left(\xix  P(\n) u_xdx+\xix  \n fudx\right)\\
  &\quad-\frac{1}{2}\xix [\mu(\n)-\mu'(\n)\n]u_x^3dx+\ga\xix P(\n)u_x^2dx-\xix \n u^2 f_xdx\\
  &\le\frac{d}{dt}\left(\xix P(\n) u_xdx+\xix \n  fudx\right)\\
  &\quad+C \|u_x\|_{L^3}^3+C \|u_x\|_{L^2}^2+C  \|u\|_{L^\infty}^2\|f_x\|_{L^2}\\
    &\le\frac{d}{dt}\left(\xix  P(\n) u_xdx+\xix \n fudx\right) +C \|u_x\|_{L^\infty}\|u_x\|_{L^2}^2+C\|u_x\|_{L^2}^2,
\ea
\ee
where one has used \eqref{gj3.2} and the following fact:
\be
\ba\la{3.5}
  \|u\|_{L^\infty}\le C\|u\|_{L^2}^{1/2}\|u_x\|_{L^2}^{1/2}\le C\|u_x\|_{L^2}.
\ea
\ee

 Then, using \eqref{1d}$_2$, \eqref{1d1}, \eqref{gj3.1}, and \eqref{gj3.2}, we have
 \be
\ba\la{3.7}
  \|u_x\|_{L^\infty} \le  &C  \|\mu(\n)u_x-P(\n)\|_{L^\infty}+C\|P(\n)\|_{L^\infty}\\
    \le &C  \|\mu(\n)u_x-P(\n)\|_{L^1}+C  \|(\mu(\n)u_x-P(\n))_x\|_{L^1}+C \\
    \le & C\|\mu(\n)u_x\|_{L^1}+C\|P(\n)\|_{L^1}+C\|\n\dot u\|_{L^1}+C\|\n f\|_{L^1}+C \\
    \le & C \xix\mu(\n)u_x^2dx+C \|\n^{1/2}\dot u\|_{L^2}+C ,
\ea
\ee
which together with \eqref{3.2} and Young's inequality gives
 \be
\ba\la{3.8}
  & \frac{d}{dt}B(t)+\frac{1}{2}\xix\n\dot u^2dx \le  C\left(\xix\mu(\n)u_x^2dx\right)^2+C\xix\mu(\n)u_x^2dx,
\ea
\ee
where
$$B(t)\triangleq   \frac{1}{2} \xix \mu(\n)u_x^2dx-\xix P(\n)u_xdx-\xix \n fudx
$$
satisfies
 \be
\ba\la{3.6}
\frac{1}{4}\|\sqrt{\mu(\n)}u_x\|_{L^2}^2-C\le  B(t)\le  C\|\sqrt{\mu(\n)}u_x\|_{L^2}^2+C
\ea
\ee
owing to the following estimate:
 \bnn\la{3.10}
   \xix P(\n) u_xdx+ \xix \n fudx
      \le \frac{1}{4}\xix \mu(\n)u_x^2dx+C,
\enn
due to \eqref{1d1}, \eqref{gj3.1}, and \eqref{gj3.2}.

Hence,  Gronwall's inequality together with \eqref{3.8}, \eqref{3.6}, and  \eqref{gj3.1} implies \eqref{gj3.5} and thus finishes the proof of Lemma \ref{la3.5}. \hfill $\Box$

\begin{lemma} \la{la3.6}
Let $(\n,u)$ be a smooth solution of
 \eqref{1d}-\eqref{n3}  on $(0,1) \times [0,T].$ Then
\be
\ba\la{gj3.6}
 \sup_{0\le t\le T}\si \left(\|u\|^2_{W^{1,\infty} }+\|\n^{1/2}\dot u\|_{L^2}^2\right) +\xiT \si\|  \dot u_x\|_{L^2}^2 dt\le C,
 \ea
\ee where and in what follows, $\si(t)\triangleq\min \{1,t\}.$
\end{lemma}

\pf First, operating $\pa_t+(u\cdot~)_x$ to $ (\ref{1d})_2 $ yields that
  \be
\ba\la{3.12}
  \n \dot u_t+\n u\dot u_x-[\mu(\n) \dot u_x]_x=-\gamma [P(\n)u_x]_x-[(\mu(\n)+\mu'(\n)\n)u_x^2]_x+\n uf_x,
\ea
\ee which
multiplied by $\dot{u} $ gives
 \be
\ba\la{3.13}
  &\frac{1}{2}\frac{d}{dt}\xix \n|\dot u|^2dx+\xix \mu(\n) |\dot u_x|^2dx\\
  &=\gamma \xix P(\n)u_x\dot u_xdx+\xix (\mu(\n)+\mu'(\n)\n)u_x^2 \dot u_xdx +\xix \n uf_x \dot udx\\
  &\le C \|u_x\|_{L^2} \|\dot u_x\|_{L^2}(1+\|u_x\|_{L^\infty} )+C \|u\|_{L^\infty}\|f_x\|_{L^2}\|\n^{1/2}\dot u\|_{L^2}\\
  &\le C \|u_x\|_{L^2} \|\dot u_x\|_{L^2}(1+ \|\n^{1/2}\dot u\|_{L^2} )+C \|u_x\|_{L^2}\|\n^{1/2}\dot u\|_{L^2}\\
  &\le \ve\|\dot u_x\|_{L^2}^2+C(\ve) \|u_x\|_{L^2}^2  +C  \|\n^{1/2}\dot u\|_{L^2}^2
\ea
\ee
due to \eqref{3.7} and  \eqref{gj3.5}.

Then, multiplying \eqref{3.13} by $\si$ and integrating the resultant inequality  over $(0,T)$, we obtain after  choosing $\ve$ suitably small and using \eqref{gj3.5}  and \eqref{gj3.1} that
 \be
\ba\la{3.16'}
   &\sup_{0\le t\le T} \si \xix \n|\dot u|^2dx+\xiT \si\xix \mu(\n) |\dot u_x|^2dxdt \le C,
\ea
\ee
which together with \eqref{3.5}, \eqref{3.7},   and \eqref{gj3.5} leads to
 \bnn\la{gj3}
   &\sup_{0\le t\le T} \si \|u_x\|^2_{L^\infty} \le C.
\enn
Combining this  with \eqref{3.16'} gives \eqref{gj3.6}. The proof of Lemma \ref{la3.6} is finished. \hfill $\Box$

 \subsection{\la{se3'} A priori estimates(II): higher order estimates}
 In this subsection, we prove the higher-order estimates of the smooth solution  $(\n,u)$ to the problem \eqref{1d}-\eqref{n3}.

 \begin{lemma} \la{la4.1} For any given $T>0$, there exists a   positive constant $C$ depending  on $T$,
 $\gamma, \bar\mu, \|\n_0\|_{H^1}, \|u_0\|_{H^1}$, and $ \|f\|_{H^1}$ such that
\be
\ba\la{gj4.1}
   \sup_{0\le t\le T} \left(\|\n_x\|_{L^2}^2+\|\n_t\|_{L^2}^2\right) \le C,
\ea
\ee
and
\be
\ba\la{ngj4.1}
   \sup_{0\le t\le T} \left(\si\|\n^{1/2}u_t\|_{L^2}^2 +\si\| u_{xx}\|_{L^2}^2\right) +\int_0^T \left( \| u_{xx}\|_{L^2}^2+\si\| u_{xt}\|_{L^2}^2\right)dt\le C.
\ea
\ee
\end{lemma}

\pf First, differentiating $\eqref{1d}_1$ with respect to $x$ gives
\be
\ba\la{a1}
   \n_{xt}+\n_{xx}u+2\n_xu_x+\n u_{xx}=0.
\ea
\ee
Multiplying \eqref{a1} by $\n_x$ and integrating the resulting equation by parts yield that
\be\ba\la{a2}
   \frac{d}{dt}\|\n_x\|_{L^2} \le C\|u_x\|_{L^\infty}\|\n_x\|_{L^2} +C \|u_{xx}\|_{L^2}
\ea
\ee
due to \eqref{gj3.2}.

Next, it is easy to deduce from $\eqref{1d}_2$ that
\be
\ba\la{a7}
\mu(\n)u_{xx}= \n \dot u+P_x-\n f-\mu'(\n)\n_xu_x,
\ea
\ee
which together with  \eqref{n3}, \eqref{gj3.2}, \eqref{gj3.5}, and \eqref{gj3.6} yields
\be
\ba\la{a9}
\|u_{xx}\|_{L^2}&\le  C\| \n^{1/2} \dot u\|_{L^2}+C\| \n_x\|_{L^2}(\| u_x\|_{L^\infty}+1)+C\\&\le C\si^{-1/2}(1+ \| \n_x\|_{L^2}).
\ea
\ee

Submitting \eqref{a9} into \eqref{a2} and using Gronwall's inequality, one gets
\be
\ba\la{a5}
   \sup_{0\le t\le T}  \|\n_x\|_{L^2}  \le C,
\ea
\ee
which along with $\eqref{1d}_1$ and \eqref{gj3.5} leads to
\be
\ba\la{a6}
 \|\n_t\|_{L^2}  & \le \|u \|_{L^\infty} \|\n_x\|_{L^2} +C\| u_x\|_{L^2} \le C.
\ea
\ee
The combination of \eqref{a5} with  \eqref{a6} leads to \eqref{gj4.1}.

Now, it follows from \eqref{gj3.2}, \eqref{3.7}, \eqref{gj3.5}, \eqref{a9}, and \eqref{gj4.1} that
\be\ba\la{b5}
 \|\n^{1/2}u_t\|_{L^2}^2&\le \|\n^{1/2}\dot u\|_{L^2}^2+ \|\n^{1/2}u u_x\|_{L^2}^2\\
 &\le \|\n^{1/2}\dot u\|_{L^2}^2+  \bar \n \|u\|_{L^\infty}^2\| u_x\|_{L^2}^2\le \|\n^{1/2}\dot u\|_{L^2}^2+ C,\\
 \|u_{xt}\|_{L^2}^2&= \|(\dot u-uu_x)_x\|_{L^2}^2\le  \|\dot u_x\|_{L^2}^2+ \|u_x^2\|_{L^2}^2+ \|uu_{xx}\|_{L^2}^2\\
 &\le \|\dot u_x\|_{L^2}^2+ \|u_x\|_{L^\infty}^2\|u_x\|_{L^2}^2+ \|u\|_{L^\infty}^2\|u_{xx}\|_{L^2}^2\\
 &\le \|\dot u_x\|_{L^2}^2+\|\n^{1/2}\dot u\|_{L^2}^2+ C.
\ea\ee
%
Hence, \eqref{ngj4.1} is a direct consequence of \eqref{b5}, \eqref{gj3.6},  \eqref{a9}, and \eqref{gj3.5}. We complete the proof of Lemma \ref{la4.1}. \hfill $\Box$

From now on, assume that $(\n,u)$ is a smooth solution of the problem \eqref{1d}-\eqref{n3} with the smooth initial data satisfying the condition in Theorem \ref{th1} and $\mu(\cdot)\in C^2[0,\infty)$. In the following of this subsection, the general constant $C$  may depend on $T$, $\gamma, \bar\mu, \|\n_0\|_{H^2}, \|P(\n_0)\|_{H^2}, \|u_0\|_{H^2}$,   $ \|f\|_{H^1}$, and $\|g\|_{L^2}$ with $g$ as in \eqref{1d4}.

\begin{lemma} \la{nla4.3}  For any given $T>0$, there exists a   positive constant $C$ such that
\be
\ba\la{gj4.2}
   \sup_{0\le t\le T} \xl(\| u_{xx}\|_{L^2} + \| u \|_{W^{1,\infty}} + \|\n^{1/2}u_t\|_{L^2}^2 \xr)+\xiT \| u_{xt}\|_{L^2}^2 dt\le C.
\ea
\ee
\end{lemma}

\pf It follows from  the compatibility condition \eqref{1d4} that we can set 
\bnn\label{gy1}\n^{1/2}\dot u(x,t=0) =g+\n_0^{1/2}f\in L^2.\enn
Integrating \eqref{3.13} over $(0,T)$, we obtain after  using \eqref{gj3.5}  and \eqref{gj3.1} that
 \bnn
   &\sup_{0\le t\le T}  \xix \n|\dot u|^2dx+\xiT  \xix \mu(\n) |\dot u_x|^2dxdt \le C,
\enn
which together with \eqref{b5}, \eqref{3.5}, \eqref{3.7},   and \eqref{gj3.5}  gives \eqref{gj4.2} and completes the proof of Lemma \ref{nla4.3}.  \hfill $\Box$

The following higher order estimates of the solutions  are used to guarantee
the extension of local classical solution to be a global one, whose proof  are similar to those in \cite{wen2011}, see also \cite{hlx1,lx1,lisiam} considering the high dimensional case. And, we also sketch them here for completeness.

\begin{lemma} \la{la4.3}  For any given $T>0$, there exists a   positive constant $C$ such that
\be
\ba\la{gj4.3}
   &\sup_{0\le t\le T} \left( \|\n_{xx}\|_{L^2}^2 + \|P_{xx}\|_{L^2}^2 +\|\n_{xt}\|_{L^2}^2 + \|P_{xt}\|_{L^2}^2 \right) \\
   &\quad\qquad+\xiT \left(\|u_{xxx}\|_{L^2}^2 + \|\n_{tt}\|_{L^2}^2 + \|P_{tt}\|_{L^2}^2 \right) dt\le C.
\ea
\ee
\end{lemma}

\pf
Differentiating $\eqref{a1}$ with respect to $x$ gives
\be
\ba\la{c1}
   (\n_{xx})_t+\n_{xxx}u+3\n_{xx}u_x+3\n_x u_{xx}+\n u_{xxx}=0.
\ea
\ee
Multiplying \eqref{c1} by $\n_{xx}$ and integrating the resulting equation  by parts, it holds that
\be
\ba\la{c2}
     \frac{1}{2}\frac{d}{dt}\|\n_{xx}\|_{L^2}^2&=-\frac{5}{2}\xix \n_{xx}^2u_xdx-3\xix \n_x\n_{xx}u_{xx}dx-\xix \n \n_{xx}u_{xxx}dx\\
     &\le C\|u_x\|_{L^\infty} \|\n_{xx}\|_{L^2}^2+C\|\n_x\|_{L^\infty}\| u_{xx}\|_{L^2}\|\n_{xx}\|_{L^2}+C\| u_{xxx}\|_{L^2}\|\n_{xx}\|_{L^2}\\
     &\le C\|\n_{xx}\|_{L^2}^2+C\| u_{xxx}\|_{L^2}^2+C,
\ea
\ee
where in the last inequality one has  used  \eqref{gj4.2}, \eqref{gj4.1}, and the following estimate:
\be
\ba\la{c3}
     \|\n_{x}\|_{L^\infty}\le  \|\n_{x}\|_{L^1}+ \| \n_{xx}\|_{L^1} \le  C + C\| \n_{xx}\|_{L^2}.
\ea
\ee

Since  $P(\n)$ satisfies
\be
\ba\la{c4}
P_t+P_xu+\gamma Pu_x=0,
\ea
\ee
following  the same arguments as \eqref{c2}, one has
\be
\ba\la{c5}
     \frac{1}{2}\frac{d}{dt}\|P_{xx}\|_{L^2}^2\le C  \|P_{xx}\|_{L^2}^2+C \| u_{xxx}\|_{L^2}^2+C.
\ea
\ee
Hence, the combination of \eqref{c5} with  \eqref{c2} implies that
\be
\ba\la{c6}
     \frac{1}{2}\frac{d}{dt}\left(\|\n_{xx}\|_{L^2}^2+\|P_{xx}\|_{L^2}^2\right)\le C\left(\|\n_{xx}\|_{L^2}^2+\|P_{xx}\|_{L^2}^2\right)+C+C\| u_{xxx}\|_{L^2}^2.
\ea
\ee

  In order to estimate $\| u_{xxx}\|_{L^2}^2$, differentiating $\eqref{a7}$ with respect to $x$ gives
\be
\ba\la{c7}
 \mu(\n)u_{xxx}= &\n_x u_t+\n_x uu_x+ \n u_{xt}+\n u_x^2+\n uu_{xx}+P_{xx}\\
 &-\n_x f-\n f_x-\mu''(\n)\n_x^2u_x-\mu'(\n)\n_{xx}u_x-2\mu'(\n)\n_xu_{xx}.
\ea
\ee
 Similar to \eqref{3.5}, we also have
\be
\ba\la{c15}
 \| u_t\|_{L^\infty}\le C \|u_{xt}\|_{L^2},
\ea 
\ee which  along with \eqref{c7}, \eqref{n3}, \eqref{gj3.2}, \eqref{gj4.1}, \eqref{c3}, and \eqref{gj4.2} implies that
\be
\ba\la{c8}
 \| u_{xxx}\|_{L^2}^2
 &\le C\|\n_x\|_{L^2}^2\|u_t\|_{L^\infty}^2+C\|\n_x\|_{L^2}^2\|u\|_{L^\infty}^2\|u_x\|_{L^\infty}^2
 +C\|u_{xt}\|_{L^2}^2+C\|u_x\|_{L^4}^4\\
 &\quad +C\|u\|_{L^\infty}^2\|u_{xx}\|_{L^2}^2+C\|P_{xx}\|_{L^2}^2+C\|\n_x\|_{L^2}^2\|f\|_{L^\infty}^2+C\|f_x\|_{L^2}^2 \\
 &\quad +C\|\n_x\|_{L^\infty}^2\|\n_x\|_{L^2}^2\|u_x\|_{L^\infty}^2+C\|u_x\|_{L^\infty}^2\|\n_{xx}\|_{L^2}^2+C\|\n_x\|_{L^\infty}^2\|u_{xx}\|_{L^2}^2\\
 &\le C \left(1+\|u_{xt}\|_{L^2}^2+\|\n_{xx}\|_{L^2}^2+\|P_{xx}\|_{L^2}^2\right).
\ea
\ee

Submitting \eqref{c8} into \eqref{c6}, one obtains after using   Gronwall's inequality and \eqref{gj4.2} that
\be
\ba\la{c11}
 \sup_{0\le t\le T} \left( \|\n_{xx}\|_{L^2}^2 + \|P_{xx}\|_{L^2}^2
\right)\le C.\ea
\ee
This together with \eqref{a1}, \eqref{gj4.2},  and \eqref{gj4.1} implies that
\be
\ba\la{c12}
 \|\n_{xt}\|_{L^2}^2\le \|\n_{xx}\|_{L^2}^2\| u\|_{L^\infty}^2+ C\|\n_{x}\|_{L^2}^2\| u_x\|_{L^\infty}^2+ C\|u_{xx}\|_{L^2}^2\le C.
\ea
\ee

Differentiating $\eqref{1d}_1$ with respect to $t$ gives
\bnn
   \n_{tt}=-\n_{xt}u-\n_xu_t-\n_tu_x-\n u_{xt}
\enn
which combined with  \eqref{c12}, \eqref{gj4.1},  \eqref{c15}, and \eqref{gj4.2}  implies that
\be
\ba\la{c14}
\xiT \|\n_{tt}\|_{L^2}^2dt&\le \xiT \left(\|\n_{xt}\|_{L^2}^2\| u\|_{L^\infty}^2+ \|\n_{x}\|_{L^2}^2\| u_t\|_{L^\infty}^2 \right)dt\\& \quad+C\xiT \left(    \|\n_{t}\|_{L^2}^2\| u_x\|_{L^\infty}^2+ C\|u_{xt}\|_{L^2}^2\right)dt\\
&\le C\xiT \left(1+  \|u_{xt}\|_{L^2}^2\right)dt\le C.
\ea
\ee
Similarly, one can get
 \bnn
 \|P_{xt}\|_{L^2}^2+\xiT \|P_{tt}\|_{L^2}^2dt\le C,
\enn
which together with  \eqref{c8}, \eqref{gj4.2},  \eqref{c11}, \eqref{c12},  and \eqref{c14} implies  \eqref{gj4.3}. The proof of Lemma \ref{la4.3} is completed. \hfill $\Box$

\begin{lemma} \la{la4.4}  For any given $T>0$, there exists a   positive constant $C $ such that
\be
\ba\la{gj4.4}
   &\sup_{0\le t\le T} \left( t\|u_{xt}\|_{L^2}^2 + t\| u_{xxx}\|_{L^2}^2 \right) +\xiT  \xl(t\|\n^{1/2}u_{tt}\|_{L^2}^2+t\|u_{xxt}\|_{L^2}^2\xr) dt\le C.
\ea
\ee
\end{lemma}

\pf Differentiating $\eqref{1d}_2$ with respect to $t$ gives
\be
\ba\la{b1}
\n u_{tt}+\n uu_{xt}-[\mu(\n)u_x]_{xt}=-\n_t(u_t+uu_x)-\n u_tu_x-P_{xt}+\n_tf.
\ea
\ee
Multiplying \eqref{b1}  by $u_{tt}$ and integrating the resulting equation   by parts lead to
\be
\ba\la{e1}
 &\frac{1}{2}\frac{d}{dt}\xix \mu u_{xt}^2dx+\|\n^{1/2}u_{tt}\|_{L^2}^2\\
 &=\frac{1}{2}\xix \mu_t u_{xt}^2dx-\xix \mu_tu_xu_{xtt}dx+\xix\n_t fu_{tt}dx-\xix\n_t\dot uu_{tt}dx\\
 &\quad-\xix \n u_tu_xu_{tt}dx-\xix \n uu_{xt}u_{tt}dx-\xix P_{xt}u_{tt}dx\\
 &\triangleq \sum_{i=1}^7 J_i.
\ea
\ee

The term $J_i$ on the right hand of \eqref{e1} can be estimated as follows.

Using   \eqref{gj4.1} and \eqref{gj4.3}, it holds
\be
\ba\la{e2}
\|\n_t\|_{L^\infty}
\le C\| \n_t\|_{L^2}+C\|\n_{xt}\|_{L^2}\le C,
\ea
\ee
which implies that
\be
\ba\la{e3}
J_1\le C\|\mu_t\|_{L^\infty}\|u_{xt}\|_{L^2}^2\le C\|\n_t\|_{L^\infty}\|u_{xt}\|_{L^2}^2 \le C\|u_{xt}\|_{L^2}^2.
\ea
\ee

It deduces from   \eqref{gj3.2}, \eqref{gj4.1}, and \eqref{e2} that
\bnn
 \|\mu_{tt}\|_{L^2}^2 &\le  \|\mu''(\n)\n_t^2\|_{L^2}^2+ \|\mu'(\n)\n_{tt}\|_{L^2}^2\\
 &\le C\|\n_t\|_{L^\infty}^2\|\n_t\|_{L^2}^2+ C\|\n_{tt}\|_{L^2}^2\\
 &\le C + C\|\n_{tt}\|_{L^2}^2,
\enn
which along with  \eqref{gj4.2} and  \eqref{e2} gives
\be
\ba\la{e4}
J_2&=-\frac{d}{dt}\xix \mu_tu_x u_{xt} dx+\xix \mu_{tt}u_x u_{xt} dx+\xix \mu_t u_{xt}^2 dx\\
&\le-\frac{d}{dt}\xix \mu_tu_x u_{xt} dx+ C\|u_x\|_{L^\infty}\|\mu_{tt}\|_{L^2}\|u_{xt}\|_{L^2} + C\|\mu_t\|_{L^\infty}\|u_{xt}\|_{L^2}^2\\
&\le  -\frac{d}{dt}\xix \mu_tu_x u_{xt} dx+C\|u_{xt}\|_{L^2}^2+C+C\|\n_{tt}\|_{L^2}^2.
\ea
\ee

Next, it is easy to derive  from \eqref{c15} that
\be
\ba\la{e6}
J_3&=\frac{d}{dt}\xix \n_t f u_{t} dx-\xix \n_{tt}fu_t dx\\
&\le \frac{d}{dt}\xix \n_t f u_{t} dx+ C\|u_t\|_{L^\infty}\|\n_{tt}\|_{L^2}\|f\|_{L^2} \\
&\le \frac{d}{dt}\xix \n_t f u_{t} dx+ C\|u_{xt}\|_{L^2}^2+C\|\n_{tt}\|_{L^2}^2.
\ea
\ee

For the term  $J_4$, we have
\be
\ba\la{e7}
J_4&=\xix (\n u)_x\dot uu_{tt}dx=-\xix \n u\dot uu_{xtt}dx-\xix \n u\dot u_xu_{tt}dx\\
&=-\frac{d}{dt}\xix \n u \dot uu_{xt}dx+\xix (\n u)_t\dot uu_{xt}dx+\xix \n u\dot u_tu_{xt}dx-\xix \n u\dot u_xu_{tt}dx\\
&\triangleq -\frac{d}{dt}\xix \n u \dot uu_{xt}dx+J_{4,1}+J_{4,2}+J_{4,3}.
\ea
\ee
One deduces from \eqref{gj4.2},  \eqref{gj4.1},   and \eqref{c15} that
\be
\ba\la{e8}\notag
 J_{4,1}&=\xix \n_t u \dot uu_{xt}dx+\xix \n u_t\dot uu_{xt}dx\\
 &\le C\xl(\|u_t\|_{L^\infty}+\|uu_x\|_{L^\infty}\xr) \|u_{xt}\|_{L^2}\|u\|_{L^\infty}\|\n_t\|_{L^2}+ C\|\n^{1/2}\dot u\|_{L^2}\|u_t\|_{L^\infty}\|u_{xt}\|_{L^2}\\
  &\le C +C\|u_{xt}\|_{L^2}^2,\\
 J_{4,2}&=\xix \n u u_{tt}u_{xt}dx+\xix \n u u_tu_xu_{xt}dx+\xix \n u^2u_{xt}^2dx\\
 &\le C\|\n^{1/2}u\|_{L^\infty} \|u_{xt}\|_{L^2} \|\n^{1/2} u_{tt}\|_{L^2}\\
 &\quad+ C\|\n u\|_{L^\infty}\|u_t\|_{L^\infty}\|u_{xt}\|_{L^2}\|u_x\|_{L^2}+C\|\n u^2\|_{L^\infty} \|u_{xt}\|_{L^2}^2\\
 &\le \ve \|\n^{1/2} u_{tt}\|_{L^2}^2+C\|u_{xt}\|_{L^2}^2,\\
 J_{4,3}&\le C\|\n^{1/2}u\|_{L^\infty} \|\dot u_{x}\|_{L^2} \|\n^{1/2} u_{tt}\|_{L^2}\le \ve \|\n^{1/2} u_{tt}\|_{L^2}^2+C\|u_{xt}\|_{L^2}^2+C.
\ea
\ee
This combined with \eqref{e7} yields
\be
\ba\la{e9}
J_4 \le -\frac{d}{dt}\xix \n u \dot uu_{xt}dx+\ve \|\n^{1/2} u_{tt}\|_{L^2}^2+C\|u_{xt}\|_{L^2}^2+C .
\ea
\ee

Moreover, by virtue of  \eqref{gj4.2}, it holds that
\be
\ba\la{e10}
J_5+J_6&\le  C\|u_x\|_{L^\infty} \|\n^{1/2}u_{tt}\|_{L^2}\|\n^{1/2}u_t\|_{L^2}+C\|\n^{1/2}u\|_{L^\infty} \|\n^{1/2}u_{tt}\|_{L^2}\|u_{xt}\|_{L^2}\\
&\le \ve \|\n^{1/2} u_{tt}\|_{L^2}^2+C\|u_{xt}\|_{L^2}^2+C.
\ea
\ee

Finally, we can estimate $J_7$ as follows,
\be
\ba\la{e11}
J_7 &=\xix P_tu_{xtt}dx= \frac{d}{dt}\xix P_tu_{xt}dx-\xix P_{tt}u_{xt}dx\\
&\le \frac{d}{dt}\xix P_tu_{xt}dx +C\|P_{tt}\|_{L^2}^2+C\|u_{xt}\|_{L^2}^2.
\ea
\ee

Submitting  \eqref{e3}-\eqref{e11} into \eqref{e1} and choosing $\ve$ suitably small, we have
\be
\ba\la{e13}
\frac{d}{dt}\Pi(t)+\|\n^{1/2}u_{tt}\|_{L^2}^2\le C\|u_{xt}\|_{L^2}^2+C +C\|\n_{tt}\|_{L^2}^2+C\|P_{tt}\|_{L^2}^2.
\ea
\ee
where
\bnn
\ba\la{e12}
\Pi(t)&\triangleq  \xix \mu u_{xt}^2dx+\xix \mu_tu_x u_{xt} dx+\xix \n u \dot uu_{xt}dx\\
&\quad-\xix P_tu_{xt}dx-\xix \n_t f u_{t} dx
\ea
\enn
satisfies
\be
\ba\la{e15}
 \frac{\bar\mu}{2}\|u_{xt}\|_{L^2}^2-C \le \Pi(t) \le C\|u_{xt}\|_{L^2}^2+C
\ea
\ee
owing to  the following estimates:
\bnn
\ba\la{e14}
 &\left|\xix \mu_tu_x u_{xt} dx+\xix \n u \dot uu_{xt}dx-\xix P_tu_{xt}dx-\xix \n_t f u_{t} dx\right|\\
&\le C\|\n_t\|_{L^2} \| u_{x}\|_{L^\infty} \|  u_{xt}\|_{L^2}+C\|\n^{1/2}u\|_{L^\infty} \|\n^{1/2}\dot u\|_{L^2} \| u_{xt}\|_{L^2}\\
&\quad+C \|P_t\|_{L^2} \| u_{xt}\|_{L^2}+C\|u_t\|_{L^\infty} \|\n_t\|_{L^2} \|f\|_{L^2}\\
&\le  \frac{\bar\mu}{2}\|u_{xt}\|_{L^2}^2+C,
\ea
\enn
where one has used \eqref{gj4.2},  \eqref{gj4.1}, and \eqref{c15}.

Hence, the Gronwall's inequality   together with \eqref{e13}, \eqref{e15}, \eqref{gj4.2}, and \eqref{gj4.3}  gives
\be
\ba\la{e16}
   \sup_{0\le t\le T}  t\|u_{xt}\|_{L^2}^2 +\xiT t\|\n^{1/2}u_{tt}\|_{L^2}^2 dt\le C,
\ea
\ee
which along with  \eqref{c8} and  \eqref{gj4.3} implies
 \be
\ba\la{e18}
    \sup_{0\le t\le T}  t\|u_{xxx}\|_{L^2}^2 \le C.
\ea
\ee

Next, it follows from \eqref{b1} that
 \bnn
\mu(\n)u_{xxt}&=\n u_{tt}+\n_t\dot u+\n u_tu_x+\n uu_{xt}+P_{xt}-\n_tf\\
&\quad-\mu''(\n)\n_t\n_xu_{x}-\mu'(\n)\n_{xt}u_{x}-\mu'(\n)\n_{x}u_{xt}-\mu'(\n)\n_{t}u_{xx}.
\enn
Combining this   with  \eqref{c15},  \eqref{gj4.1}, \eqref{gj4.2}, \eqref{gj4.3}, and \eqref{e2} yields that
 \be
\ba\la{Re20}\notag
\|u_{xxt}\|_{L^2}^2 
&\le C\|\n^{1/2} u_{tt}\|_{L^2}^2+ C\|\dot u\|_{L^\infty}^2\|\n_t\|_{L^2}^2+C\|u_x\|_{L^\infty}^2\|\n^{1/2}u_t\|_{L^2}^2\\
&\quad+C\|\n u\|_{L^\infty}^2\|u_{xt}\|_{L^2}^2+C\|P_{xt}\|_{L^2}^2+C\|f\|_{L^\infty}^2\|\n_t\|_{L^2}^2+C\|\n_t\|_{L^\infty}^2\|u_{xx}\|_{L^2}^2\\
&\quad+C\|u_x\|_{L^\infty}^2\|\n_t\|_{L^\infty}^2\|\n_x\|_{L^2}^2+C\|u_x\|_{L^\infty}^2\|\n_{xt}\|_{L^2}^2+C\|\n_x\|_{L^\infty}^2\|u_{xt}\|_{L^2}^2\\
&\le C\|\n^{1/2} u_{tt}\|_{L^2}^2+ C\|u_{xt}\|_{L^2}^2+ C,
\ea
\ee
which along with  \eqref{e16} implies that
 \be
\ba\la{e21}
\xiT t\|u_{xxt}\|_{L^2}^2dt\le  C.
\ea
\ee
The combination of  \eqref{e16}--\eqref{e18}  with  \eqref{e21} yields \eqref{gj4.4} and   completes   the proof of Lemma \ref{la4.4}. \hfill $\Box$

\section{Proof of Theorems \ref{th11} and \ref{th1}}

{\it Proof of Theorem \ref{th11}.}  With all the a priori estimates obtained in section 2 at hand, we will divide the proof into three steps.

\textbf{Step 1}. We prove the local existence and uniqueness of the strong solution when the initial density contains vacuum. That is, Theorem \ref{th11} holds for some $T_0>0$.

Let  $(\n_{0},u_{0},f)$  be as  in Theorem \ref{th11},  we construct
\be\label{gy2}\n_0^{\delta}=\hat\n_0^\de+\delta,~~u_0^{\delta}=u_0\ast j_\delta,~~f^{\delta}=f\ast j_\delta,\ee
where $j_\delta$ is the standard mollifying kernel of width $\delta$ and  $0\le \hat\n_0^\de\in  C_0^\infty (0,1)$ satisfies
   \be\label{ma5.2} \ba \hat\n_0^\de\rightarrow  \n_{0} ~\quad {\rm in}\,\, H^1,~~~  {\rm as}~\de\rightarrow 0.\ea\ee
   Thus, we have
      \be\label{ma6.1} \ba  &\n_0^{\delta}\rightarrow  \n_{0},~~u_0^{\delta}\rightarrow  u_{0},~~f^\delta\rightarrow  f, ~~&\quad {\rm in}\,\, H^1, ~~~  {\rm as}~\de\rightarrow 0,\ea\ee
      and
\be\label{gyma6.1} \ba \|\n_0^{\delta}\|_{H^1} \le  C+C\|\n_0 \|_{H^1},  \,\,  \|u_0^{\delta}\|_{H^1}\le C\|u_0\|_{H^1}, \,\,  \|f^{\delta}\|_{H^1}\le C\|f\|_{H^1}.\ea\ee

 By virtue of  Lemma \ref{le3}, the initial boundary problem \eqref{1d}-\eqref{n3}  with the initial data $(\n_0^{\delta},  u_0^\delta)$ has a classical solution  $(\n^{\delta},  u^\delta)$ on $(0,1)\times [0,T_0]$. Furthermore, the estimates obtained  in  Lemmas \ref{la3.2}--\ref{la4.1}   show that the solution $(\n^{\delta},  u^\delta)$ satisfies for any $0<T<+\infty$,
\be\ba
   \la{8.6}
&  \sup_{0\le t\le T} \xl(\|(\n^\delta,\mu(\n^\delta),P(\n^\delta))\|_{H^1}+\| \n^\delta_t\|_{L^2}+\| \n^\delta u^\delta \|_{L^2}+ \|u^\delta\|_{H^1}\right.\\
&\quad\quad\quad\left.+\sqrt{t}\|\sqrt{\n^\delta}u^\delta_t\|_{L^2}+\sqrt{t}\|u^\delta_{xx}\|_{L^2}\xr) +\xiT\xl( \|u^\delta\|_{H^2}^2+ t\|u^\delta_{xt}\|_{L^2}^2\xr)dt\le \bar C, \ea\ee
where $\bar C$ is independent of $\delta$.   With all the estimate \eqref{8.6}  at hand,  we  find  that the sequence
$(\n^{\delta},u^{\delta})$ converges, up to the extraction of subsequences, to some limit $(\n,u)$   in the
obvious weak sense.  Then letting $\delta\rightarrow 0$, we deduce from \eqref{8.6} that  $(\n,u)$ is a strong solution of \eqref{1d}-\eqref{n3}   on $(0,1)\times (0,T_0]$ satisfying
\be
   \la{all}\begin{cases}
 \n \in L^{\infty}(0,T_0;\,H^1),~~\n_t \in L^{\infty}(0,T_0;\,L^2),\\
    u\in L^{\infty}(0,T_0; \, H^1_0)\cap L^2(0,T_0; \,H^2) \\
 t^{1/2}u\in L^{\infty}(0,T_0; \,H^2),~~ t^{1/2}u_t\in L^2(0,T_0; \,H^1),\\
  t^{1/2}\n^{1/2}u_{t}\in L^{2}(0,T_0;\,L^2).\end{cases} \ee

 Then, the uniqueness of the strong solution $(\n,u)$ is guaranteed by the regularities \eqref{all}. For the detailed proof,  please see \cite{coi1,liliang}.

  \textbf{Step 2}.  We will extend the local existence time $T_0$ of  strong solution to be infinity and thus prove the global existence result.

Let $T^*$ be the  maximal time of existence for the strong solution. Then, $T^*\ge T_0$. For any $0<\tau< T\le T^*$ with $T$ finite, one deduce from
$$ u\in L^{\infty}(0,T; \, H^1_0)\cap L^2(0,T; \,H^2),~~u_t\in L^2(\tau,T;\,H^1),$$
that
 \be \la{x3} u \in C([\tau ,T]; \,  H^1_0).\ee
Furthermore, it follows from
$$\n \in L^{\infty}(0,T;\,H^1),~~\n_t \in L^{\infty}(0,T;\,L^2)$$ that
  \be \la{nx2} \n \in C([0 ,T];H^1).\ee

 Defining
  \bnn (\n^*, u^*)\triangleq(\n,   u)(x,T^*)=\lim_{t\rightarrow T^*}(\n, u)(x,t),\enn
we   derive from \eqref{x3} and  \eqref{nx2} that $ (\n^*, u^*)$ satisfies the initial condition \eqref{1d3a} at $t=T^*$.

Hence, we take $ (\n^*, u^*)$ as the initial data at $t=T^*$ and then use the local existence theory to extend  the strong solution beyond the maximum existence time $T^*$. This contradicts the assumption on $T^*$. We finally show that $T^*$ could be  infinity and prove the global existence of the strong solution.

  \textbf{Step 3}. It remains   to prove \eqref{1.12'}. Direct calculations lead to
\bnn\ba\label{ts7}
\left| \frac{d}{dt}\|u_x\|_{L^2}^2\right|& = \left|2\xix u_xu_{xt}dx\right|\\&=\left|2\xix u_x\left(\dot u_x-(uu_x)_x\right)dx \right|\\
 & =\left|2\xix u_x \dot u_xdx- \xix u_x^3 dx\right|\\
 &\le  C\|\dot u_x\|_{L^2}^2+C(1+\|u_x\|_{L^\infty})\|u_x\|_{L^2}^2,
\ea\enn
which together with \eqref{gj3.1} and \eqref{gj3.6} yields
\bnn\ba\label{ts8}
\int_1^\infty\left(\|u_x\|_{L^2}^2+\left|\frac{d}{dt}\|u_x\|_{L^2}^2\right|\right)dt\le C.
\ea\enn
Thus,
\bnn\label{ts3}
\lim_{t\rightarrow \infty}\|u_x\|_{L^2}^2(t)=0,
\enn
which combined with \eqref{gj3.6}  gives
\be\la{gj3.4'}
\lim_{t\rightarrow \infty} \|u\|_{W^{1,p}} =0,~~~~\forall~p\in[1,\infty).
\ee

The proof of Theorem \ref{th11} is finished.
\hfill$\Box$

{\it Proof of Theorem \ref{th1}.}  With the higher-order estimates in Lemmas \ref{nla4.3}--\ref{la4.4} at hand, the proof of Theorem \ref{th1} is similar to those of Theorem \ref{th11} and is omitted here for simplicity. \hfill$\Box$

\section{Proof of Theorem  \ref{th2}}

 The proof of Theorem \ref{th2} is divided into two steps as follows.

\textbf{Step 1}. We will prove
\be\ba\la{gj3.4}   \lim_{t\rightarrow0} \|\n(\cdot,t)-\n_s(\cdot)\|_{L^p}\rightarrow 0,~~\forall~p\in[1,\infty).\ea\ee

Considering the function
\bnn\label{st3.3}
\mathcal{P}(t)=\xix[P(\rho)-\tilde{P}(\n)]^2dx
\enn
with
\bnn\label{st3.2}
 \tilde{P}(\n)= \tilde{P}(\rho(x,t))\triangleq\int_0^1 P(\rho)dx+\int_0^x\rho fdy-\int_0^1\int_0^x\rho fdydx,
\enn
we claim that
\begin{equation}\label{st3.4}
    \lim_{t\rightarrow\infty}\int_{t-1}^t\left(\mathcal{P}(\tau)+\left|\frac{d}{d\tau}\mathcal{P}(\tau)\right|\right)d\tau=0.
\end{equation}

With \eqref{st3.4} at hand, one can derive the desired \eqref{gj3.4} with the same arguments as those in \cite{is2002,is2003}. For reader's convenience, we sketch them here for completeness. Indeed,
for any $t>1$ and $s\in(t-1,t)$, it holds
\bnn \label{st15}
    \mathcal{P}(t) \leq\int_{t-1}^t\left(\mathcal{P}(s)+\left|\frac{d}{d\tau}\mathcal{P}\right|\right)d\tau,
\enn
which together with \eqref{st3.4} implies that
\be\label{st3.8}
\lim_{t\rightarrow\infty}\mathcal{P}(t)=\lim_{t\rightarrow\infty}\|P(\n(\cdot, t))-\tilde{P}(\n(\cdot, t))\|_{L^2}^2=0.
\ee
Let $t_n\rightarrow\infty$ be an arbitrary sequence, the uniform upper bound of the density \eqref{gj3.2} imlpies that there exist  a function $\widetilde{\rho}\in L^\infty$ and a positive constant $P_c$ such that for some subsequence $\{t_n'\}\subset\{t_n\}$,
\begin{equation}\label{st3.9}
    \rho(\cdot,t_n')\rightarrow\widetilde{\rho}(\cdot)\text{  weakly * in $L^\infty$},~~~~
    \int_0^1 P(\rho(\cdot,t_n'))dx\rightarrow P_c.
\end{equation}
Clearly, $0\leq\widetilde{\rho}\leq C$. The standard compactness argument together with \eqref{st3.9} yields
\begin{equation}\label{st3.10}
    \tilde{P}\left(\rho(\cdot,t_n')\right)\rightarrow P_s(\cdot)\triangleq P_c
    +\int_0^x\widetilde{\rho}fdx-\int_0^1\int_0^x\widetilde{\rho}fdydx~~~\text{  in~~$C([0,1])$},
\end{equation}
which along with \eqref{st3.8} leads to
\bnn\label{st16}
P(\rho(\cdot,t_n'))\rightarrow P_s(\cdot) ~~~\text{  in~~$L^2$}.
\enn
Consequently, it holds that for some subsequence  $\{t_n''\}\subset\{t_n'\}$,
 $$P(\rho(x,t_n''))\rightarrow P_s(x) ~~~~\text{a.e.~ in}~(0,1).$$
  The continuity and monotonicity of $P(\cdot)$ deduce that the inverse function of $P$, denoted by $P^{-1}$,  is continuous,  and then
$$\rho(x,t_n'')\rightarrow \rho_s(x)\triangleq P^{-1}(P_s(x)) ~~~~\text{a.e.~ in}~(0,1),$$
where $\rho_s\in C([0,1])$ and $\rho_s >0$.
 This together with Lebesgue dominated theorem and \eqref{gj3.2} implies that
\begin{equation}\label{st3.11}
    \rho(\cdot,t_n'')\rightarrow\rho_s(\cdot),~~\text{ in $L^p$},~~~
    \forall~p \in[1,\infty),
\end{equation}
 which along with \eqref{st3.9} yields $$\rho_s=\widetilde{\rho}.$$
 According to the definition of $P_s$ in \eqref{st3.10}, it holds
$$P(\rho_s)=P_s=P_c
    +\int_0^x\rho_sfdx-\int_0^1\int_0^x\rho_sfdydx,$$
which together  with \eqref{st3.11} and \eqref{03.2} yields that
\be\notag [P(\rho_s)]_x= \rho_sf,~~~ \xix \n_s dx=1.\ee
Hence, we show that   $\rho_s$ is indeed the solution to the stationary problem \eqref{1s} due to Lemma \ref{le2}. And, \eqref{gj3.4} is a direct consequence of  \eqref{st3.11}.

Now, it remains to prove \eqref{st3.4}.
Denoting
\bnn\label{st1}
\Psi=\int_0^x (P(\rho)-\tilde{P}(\n))dx,
\enn which satisfies  $\Psi(0,t)=\Psi(1,t)=0,$  after using integration by parts and \eqref{1d}$_2$, we rewrite $\mathcal{P}(t)$    as
\be\ba\label{st4}
\mathcal{P}(t)&=\int_0^1(P(\rho)-\tilde{P}(\n))d\Psi  =-\int_0^1(P(\rho)-\tilde{P}(\n))_x\Psi dx\\
    &=\int_0^1\left((\rho u)_t+(\rho u^2)_x-[\mu(\rho)u_x]_x\right)\Psi dx\\
    &=\frac{d}{dt}\int_0^1\rho u\Psi dx-\int_0^1\rho u\Psi_t dx-\xix\left(\rho u^2-\mu(\rho)u_x\right)\Psi_x dx.
\ea\ee

 It follows from \eqref{gj3.1} and \eqref{gj3.2} that
\be\ba\label{st2}
|\Psi|+|\Psi_x|&=\left|\int_0^x\left(P(\rho)-\int_0^1P(\rho)dy-\int_0^x\rho fdy+\int_0^1\int_0^\xi\rho fdyd\xi\right)dx\right|\\
   &\quad+\left|P(\rho)-\int_0^1P(\rho)dx-\int_0^x\rho fdy+\int_0^1\int_0^x\rho fdydx\right|\leq C,
\ea\ee
and thus\be\ba\label{st5}
 \left|-\xix\left(\rho u^2-\mu(\rho)u_x\right)\Psi_x dx \right|\le C\|u_x\|_{L^2}^2+C\|u_x\|_{L^2}.
\ea\ee

Since
$$P_t+(Pu)_x+(\gamma-1)Pu_x=0,$$
one deduces from integration by parts and \eqref{1d}$_1$ that
\bnn\ba\label{st6}
   &(P(\rho)-\tilde{P}(\n))_t\\
   &=-(Pu)_x-(\gamma-1)Pu_x-\int_0^1 P_tdx-\int_0^x\rho_t fdy+\int_0^1\int_0^x\rho_t fdydx\\
   &=-(Pu)_x-(\gamma-1)Pu_x+(\gamma-1)\int_0^1Pu_xdx+\int_0^x(\n u)_y fdy-\int_0^1\int_0^x(\n u)_y fdydx\\
   &=-(Pu)_x-(\gamma-1)Pu_x+(\gamma-1)\int_0^1Pu_xdx\\
   &\quad+\n uf-\int_0^x \n u  f_ydy-\int_0^1 \n u  f dx+\int_0^1\int_0^x \n u  f_ydydx.
\ea\enn
Combining this, \eqref{gj3.2}, and \eqref{3.5}  gives
\be\ba\label{st7}
    |\Psi_t| &=\left|\int_0^x(P(\rho)-\tilde{P}(\n))_tdx\right|\\
    &=\left|-\rho^\gamma u-(\gamma-1)\int_0^x\rho^\gamma u_xdx+(\gamma-1) \int_0^1\rho^\gamma u_xdx\right.\\
   &\quad +\left.\int_0^x\left(\n uf-\int_0^x \n u  f_ydy-\int_0^1 \n u  f dx+\int_0^1\int_0^x \n u  f_ydydx\right)dx \right|\\
   &\le C(\bar\n,\|f\|_{H^1})\|u\|_{L^2}+C(\bar\n)\|u\|_{L^\infty}+C(\bar\n)\|u_x\|_{L^2}\le C\|u_x\|_{L^2}.
\ea\ee
One thus gets
\be\ba\label{st8}
 \int_0^1 \rho u\Psi_tdx\leq C\|u_x\|_{L^2}\|u\|_{L^\infty}\xix\rho dx\leq C\|u_x\|_{L^2}^2.
\ea\ee

Hence, on the one hand, it follows from \eqref{st2},  \eqref{st4}, \eqref{st5}, and \eqref{st8} that
\be\ba\label{st9}
 \int_{t-1}^t \mathcal{P}(\tau)d\tau&=\int_0^1\rho u\Psi dx(t)-\int_0^1\rho u\Psi dx(t-1)\\
 &\quad- \int_{t-1}^t\int_0^1\rho u\Psi_t dxd\tau- \int_{t-1}^t\xix\left(\rho u^2-\mu(\rho)u_x\right)\Psi_x d\tau\\
 &\le C \|u\|_{L^2}+C\int_{t-1}^t\left(\|u_x\|_{L^2}^2+\|u_x\|_{L^2}\right)d\tau\\
 &\le C \|u_x\|_{L^2}+C\int_{t-1}^t\|u_x\|_{L^2}^2d\tau+C\left(\int_{t-1}^t\|u_x\|_{L^2}^2d\tau\right)^{1/2}.
\ea\ee
On the other hand, it deduces from \eqref{st2} and \eqref{st7} that
\be\ba\label{st11}
\int_{t-1}^t \left|\frac{d}{d\tau}\mathcal{P}(\tau)\right|d\tau
&=2\int_{t-1}^t \left|\xix \Psi_x(P(\rho)-\tilde{P}(\n))_\tau dx\right|d\tau\\
& \le C\int_{t-1}^t \|u_x\|_{L^2} d\tau\leq C\left(\int_{t-1}^t\|u_x\|_{L^2}^2d\tau\right)^{1/2}.
\ea\ee

 Then, the desired \eqref{st3.4} is deduced directly from \eqref{st9}, \eqref{st11}, \eqref{gj3.4'}, and \eqref{gj3.1}.

\textbf{Step 2}.  Now, we are  in a position to prove \eqref{ed1}. The method  used here is  motivated by Huang-Li-Xin \cite{jmpa06} and Li-Zhang-Zhao \cite{lisiam}.

Thanks to \eqref{1s}, the momentum equation \eqref{1d}$_2$ can be rewritten as
\bnn\ba\label{ed1.1}
 \rho u_t + \rho u u_x + [P(\n)]_x-\n\n_s^{-1}[P(\n_s)]_x=[\mu(\n)u_x]_x,
\ea\enn
which multiplied by $u$, we obtain after using integration by parts and  \eqref{1d}$_1$ that
\be\ba\label{ed1.2}
 \frac{d}{dt}\xl(\frac{1}{2}\xix \n u^2dx+\xix G(\n)dx\xr)+\xix \mu(\n) u_x^2dx=0,
\ea\ee
where
\be\ba\label{ed2.0}
G(\n)&\triangleq \int_{\n_s}^\n  \int_{\n_s}^r \frac{P'(\xi)}{\xi}d\xi dr =\n  \int_{\n_s}^\n \frac{P(\xi)-P(\n_s)}{\xi^2} d\xi\\
&=\frac{1}{\gamma-1}\xl(P(\n)-P(\n_s)-P'(\n_s)(\n-\n_s)\xr).
\ea\ee
Clearly,  it follows from \eqref{gj3.2} and \eqref{pre1} that there are positive constants $M_1$ and $M_2$ depending only on $\gamma$, $\bar\n$, $K_1$, and $K_2$ such that
\be\ba\label{ed2.1}
 M_1(\n-\n_s)^2\le G(\n)\le M_2(\n-\n_s)^2. \ea\ee

Next, it follows from  \eqref{1d}$_2$  and \eqref{1s} that
\be\ba\label{ed2.3}
-[P(\n)-P(\n_s)]_x+\n_s^{-1}(\n-\n_s)P'(\n_s) (\n_s)_x=(\n u)_t+(\n u^2)_x-[\mu(\n)u_x]_x.
\ea\ee
For  $\Phi (x,t)\triangleq\int_0^x(\n-\n_s)dy,$
it holds that \be\ba\label{ed2.7}\Phi (0,t)= \Phi (1,t)=0\ea\ee owing to \eqref{1s} and \eqref{03.2}.
Multiplying \eqref{ed2.3} by   $\n_s^{-1}\Phi (x,t)$ gives
\be\ba\label{ed2.5}
&-\xix \xl[\n_s^{-1}(P(\n)-P(\n_s))\xr]_x\Phi  dx\\
&=\xix \n_s^{-2}(\n_s)_x \xl[P(\n)-P(\n_s)-P'(\n_s)(\n-\n_s) \xr]\Phi  dx\\
&\quad+\xix \n_s^{-1}(\n u)_t\Phi  dx+\xix \n_s^{-1}(\n u^2)_x\Phi  dx-\xix \n_s^{-1} [\mu(\n)u_x]_x\Phi  dx\\
&\triangleq \sum_{i=1}^4\bar J_i.
\ea\ee

First, integration by parts combined with \eqref{ed2.7} and \eqref{pre1} gives
\be\ba\label{ed2.6}
-\xix \xl[\n_s^{-1}(P(\n)-P(\n_s))\xr]_x\Phi  dx 
&=\xix  \n_s^{-1}P'(\hat\n)(\n-\n_s)^2dx \\&\ge C_1 \|\n-\n_s\|_{L^2}^2,
\ea\ee
where $0<\hat\n \in (\min\{\n, \n_s\}, \max\{\n, \n_s\})$ and $C_1$ is a positive constant dependent of $\gamma$, $\bar \n$, $K_1$, and $K_2$.

Next, the terms on the right hand of \eqref{ed2.5} can be estimated as follows. On the one hand, it follows from \eqref{gj3.4} that there is some $T^*>1$ such that for $t>T^*$,
\be\ba\label{ed1.8}\notag
 \|\n-\n_s\|_{L^2}  \le \frac{C_1}{4},
\ea\ee
which along with  \eqref{pre1} and \eqref{ed2.0} yields that
\be\ba\label{ed2.8}
\bar J_1 
 \le C\|\n-\n_s\|_{L^2}\|\n-\n_s\|_{L^2}^2
&\le \frac{C_1}{4}\|\n-\n_s\|_{L^2}^2.
\ea\ee
On the other hand, the integration by parts  together with  \eqref{pre1} and \eqref{gj3.2} implies that
\be\ba\label{ed2.9}
\bar J_2 &= \frac{d}{dt}\xix  \n_s^{-1}(\n u) \int_0^x(\n-\n_s)dy dx -\xix \n_s^{-1} \n u  \int_0^x \n_t dy dx \\
&=\frac{d}{dt}\xix  \n_s^{-1}(\n u) \int_0^x(\n-\n_s)dy dx+\xix \n_s^{-1} \n^2 u^2 dx\\
& \le \frac{d}{dt}\xix  \n_s^{-1}(\n u) \int_0^x(\n-\n_s)dy dx+ C\|u_x\|_{L^2}^2,
\ea\ee
and
\be\ba\label{ed2.10}
\bar J_3 +\bar J_4 &=  \xix \n_s^{-2} (\n_s)_x \n u^2 \int_0^x(\n-\n_s)dy dx- \xix \n_s^{-1}   \n u^2 (\n-\n_s) dx\\
&\quad+\xix \n_s^{-1}  \mu(\n)u_x (\n-\n_s)dx-\xix \n_s^{-2} (\n_s)_x  \mu(\n)u_x \int_0^x(\n-\n_s)dy dx\\
&  \le C \|u_x\|_{L^2}^2+C \|u_x\|_{L^2} \|\n-\n_s\|_{L^2}\\
&\le  \frac{C_1}{4} \|\n-\n_s\|_{L^2}^2+C\|u_x\|_{L^2}^2.
\ea\ee
 Substituting \eqref{ed2.6}--\eqref{ed2.10} into \eqref{ed2.5} derives
 \be\ba\label{ed3.0}
C_1\|\n-\n_s\|_{L^2}^2 \le 2\frac{d}{dt}\xix  \n_s^{-1}(\n u) \int_0^x(\n-\n_s)dy dx+ C\|u_x\|_{L^2}^2,
\ea\ee

Since
 \be\ba\label{ed3.4}\notag
\xl|2 \xix  \n_s^{-1}(\n u) \int_0^x(\n-\n_s)dy dx\xr|
&\le   C\|\sqrt{\n} u\|_{L^2}^2+  C\|\n-\n_s\|_{L^2}^2,
\ea\ee
  adding \eqref{ed3.0} multiplied  by some suitably small $\eta$ to \eqref{ed1.2} gives
 \be\ba\label{ed3.1}
&\frac{d}{dt}W(t)+\eta C_1\|\n-\n_s\|_{L^2}^2 +\frac{1}{2}\xix \mu(\n) u_x^2dx\le 0,
\ea\ee
where
 \bnn
 W(t)\triangleq \frac{1}{2}\xix \n u^2dx+\xix G(\n)dx-2\eta \xix  \n_s^{-1}(\n u) \int_0^x(\n-\n_s)dy dx
\enn
satisfies
 \be\ba\label{ed3.5}
\frac{1}{4}\xix \n u^2dx+\frac{M_1}{2}\|\n-\n_s\|_{L^2}^2\le W(t)\le C\|\sqrt{\n} u\|_{L^2}^2+ C\|\n-\n_s\|_{L^2}^2.
\ea\ee
Furthermore, one has
 \be\ba\label{ed3.2}
 \xix \n u^2dx \le C\|u \|_{L^2}^2\le \frac{C_2^{-1}}{2}\|\sqrt{\mu(\n)}u_x \|_{L^2}^2,
\ea\ee
where $C_2>0$ is a positive constant depending on $\bar \mu$.

The combination of \eqref{ed3.1} with \eqref{ed3.2} gives
\be\ba\label{ed3.6}
&\frac{d}{dt}W(t)+\al \|\n-\n_s\|_{L^2}^2 +\al  \xix \n u^2dx \le 0,
\ea\ee
where $\al$ is a positive constant depending on $\eta$, $C_1$, and $C_2$. Hence, Gronwall's inequality combined with \eqref{ed3.5} and \eqref{ed3.6} shows
\be\ba\label{ed3.7}
\|\sqrt{\n} u\|_{L^2}^2+  \|\n-\n_s\|_{L^2}^2 \le Ce^{-\al t},~~~\mbox{for}~t>T^*.
\ea\ee

In what follows, we will prove the exponential decay rate for the  $L^2$-norm of $u_x$.
 First, multiplying \eqref{ed1.2} by $e^{\frac{\al}{2}t}$, we get after using \eqref{ed2.1} and \eqref{ed3.7} 
\be\ba\label{ed4.0}
 &\frac{d}{dt}\xl(e^{\frac{\al}{2}t}\|\sqrt{\n} u\|_{L^2}^2+2e^{\frac{\al}{2}t}\xix G(\n)dx\xr)+e^{\frac{\al}{2}t}\|u_x\|_{L^2}^2\\
 &\le Ce^{\frac{\al}{2}t}\|\sqrt{\n} u\|_{L^2}^2+Ce^{\frac{\al}{2}t}\|\n-\n_s\|_{L^2}^2\\
 &\le Ce^{-\frac{\al}{2}t}.
\ea\ee
Integrating \eqref{ed4.0} over $[T^*,t]$ gives
\be\ba\label{ed4.1}
e^{\al t}\|\sqrt{\n} u\|_{L^2}^2+e^{\al t}\|\n-\n_s\|_{L^2}^2 +\int_{T^*}^te^{\frac{\al}{2}t}\|u_x\|_{L^2}^2dt\le C,~~~\mbox{for}~t>T^*,
\ea\ee  due to \eqref{ed3.7}.

Next, multiplying \eqref{ed2.3} by $\dot u$ and integrating the resulting equation by parts, it holds for $t>T^*$,
\be\ba\label{ed4.2}
 &\frac{d}{dt}\xl(\|\sqrt{\mu(\n)} u_x\|_{L^2}^2\xr)+\|\sqrt{\n} \dot u\|_{L^2}^2\\
 &=-\frac{1}{2}\xix [\mu(\n)-\mu'(\n)\n]u_x^3dx+\xix [P(\n)-P(\n_s)] \dot u_xdx\\
 &\quad +\xix \n_s^{-1}(\n-\n_s)P'(\n_s) (\n_s)_x\dot udx \\
 &\le C\|u_x\|_{L^\infty}\|u_x\|_{L^2}^2+C\|\n-\n_s\|_{L^2}\xl(\|\dot u\|_{L^2}+\|\dot u_x\|_{L^2}\xr)\\
  &\le C\|u_x\|_{L^2}^2+C\|\n-\n_s\|_{L^2}\|\dot u_x\|_{L^2},
\ea\ee
where in the last inequality one has used \eqref{gj3.6} and Poincar\'e inequality. Integrating \eqref{ed4.2} multiplied  by $e^{\frac{\al}{4}t}$ over $[T^*,t]$, we thus obtain after using \eqref{ed4.1} and \eqref{gj3.6} that \bnn
 & e^{\frac{\al}{4}t}\|u_x\|_{L^2}^2 +\int_{T^*}^te^{\frac{\al}{4}t}\|\sqrt{\n} \dot u\|_{L^2}^2dt  \le C,~~~\mbox{for}~t>T^*,
\enn which together with    \eqref{ed4.1},   \eqref{gj3.4'},  and \eqref{gj3.4} gives  \eqref{ed1}.

Finally, the proof of \eqref{bp} is similar as that of \cite[Theorem 1.2]{lijde}(see also \cite{hlx1,jmpa06}).
The proof of  Theorem \ref{th2} is completed.  \hfill$\Box$


\begin{thebibliography} {99}




\bibitem{zl89}  A. A. Amosov,  A. A. Zlotnik,  Global generalized solutions of the equations of the one-dimensional motion of a viscous heat-conducting gas, \emph{Soviet Math. Dokl.}, \textbf{38}(1989), 1-5.

\bibitem{v1989}  H.  Beir\~ao da Veiga, Long time behavior for one-dimensional motion of a general barotropic viscous fiuid, \textit{Arch. Ration.
Mech. Anal.}, \textbf{108}(1989), 141-160.





\bi{coi1} Y. Cho, H. Kim,  On classical solutions of the compressible Navier-Stokes equations with
nonnegative initial densities. \emph{Manuscript Math.},{\bf 120} (2006), 91-129.

\bibitem{wen2011} S. J. Ding, H. Y. Wen,  C. J. Zhu. Global classical large solutions of 1D compressible Navier-Stokes equations with density-dependent viscosity and vacuum, \textit{J. Differ. Eqs.},   \textbf{221}(2011), 1696-1725.





 \bibitem{fz}  D. Y. Fang,  T. Zhang, Compressible Navier-Stokes equations with vacuum state in the case of general pressure law, \emph{Math.
Methods Appl. Sci.}, \textbf{29}(2006), 1081-1106



\bibitem{jxz} S. Jiang, Z. P  Xin, P. Zhang,  Global weak solutions to 1D compressible isentropic Navier-Stokes equations with density-dependent viscosity, \emph{Methods Appl. Anal.}, \textbf{12}(2005), 239-251.







\bibitem{Kaz} Y. I. Kanel,   On a model system of equations of one-dimensional gas
motion. Differential Equations, 4 (1968),   374-380.
 \bibitem{kawohl} B. Kawohl,  Global existence of large solutions to initial boundary value problems for a viscous, heat-conducting, one-dimensional real gas, \textit{J. Differ. Eqs.},   \textbf{58}(1985), 76-103.

\bibitem{Hof} D. Hoff,  Global existence for 1D, compressible, isentropic Navier-Stokes equations with large initial data. \textit{Trans. Amer. Math. Soc.} \textbf{303} (1987), no. 1, 169--181.

\bibitem{H3}D. Hoff,  Global solutions of the Navier-Stokes equations
 for multidimensional compressible flow with discontinuous initial data.
\textit{J. Differ. Eqs.}  \textbf{120}(1995), no. 1, 215--254.

\bibitem{jmpa06} F. M. Huang, J. Li, Z.P. Xin,   Convergence to equilibria and blowup behavior of global strong solutions to the Stokes approximation equations for two-dimensional compressible flows with large data. \emph{J. Math. Pures Appl.} \textbf{(9)8}(2006), 471--491.







 \bibitem{hlx1}
  X. D. Huang; J. Li; Z. P. Xin,  Global well-posedness of classical solutions with large
 oscillations and vacuum to the three-dimensional isentropic
 compressible Navier-Stokes equations.  {Comm. Pure Appl. Math.} \textbf{65}(2012), 549-585.


\bibitem{ljcmp2008}  H. L. Li, J. Li,  Z. P. Xin,  Vanishing of vacuum states and blow-up phenomena of the compressible Navier-Stokes equations, \emph{Comm. Math. Phys.}, \textbf{281}(2008), 401-444.

\bibitem{lijde} J. Li,   Z. P. Xin, Some unifrom estimates and blowup behavior of global strong solution
to the  Stokes approximation equations for two-dimensional compressible flows, \emph{J. Differential Equations}, 221(2006), 275-308.




\bibitem{lx1} J. Li, Z. P. Xin,  Global well-posedness and large time asymptotic behavior of classical solution to the compressible Navier-Stokes equations with vacuum, http://arxiv.org/abs/1310.1673v1.

\bibitem{liliang} J. Li,  Z. L. Liang,  On classical solutions to the Cauchy problem of the two-dimensional barotropic compressible Navier-Stokes equations with vacuum. \emph{J. Math. Pures Appl.}, \textbf{102}(2014), 640-671.

\bibitem{lij1d} J. Li,  Z. L. Liang,  Some uniform estimates and large-time behavior for one-dimensional compressible Navier-Stokes system in unbounded domains with large data.  \emph{Arch. Rational Mech. Anal}. \textbf{220}(2016),  1195-1208.

  \bibitem{lisiam}  J. Li,  J. W. Zhang, J. N. Zhao,   On the global motion of viscous compressible barotropic flows subject to large external potential forces and vacuum. \emph{SIAM J. Math. Anal.} \textbf{47}(2015), no. 2, 1121--1153.

\bibitem{lxy} T. P. Liu, Z. P. Xin, T. Yang, Vacuum states of compressible flow, \textit{Discrete Contin. Dyn. Syst.}, \textbf{4}(1998), 1-32.




\bibitem{ok} M. Okada, S. Matusu-Necasov\'a, T. Makino, Free boundary problem for the equation of one-dimensional motion of compressible gas with density-dependent viscosity, \emph{Ann. Univ. Ferrara Sez.}, VII (N.S.), \textbf{48}(2002), 1-20.

\bibitem{qyz} X. Qin, Z. A. Yao, H. Zhao, One dimensional compressible Navier-Stokes equations with density-dependent viscosity and free
boundary, \emph{Commun. Pure Appl. Anal.}, \textbf{7}(2008), 373-381.





\bibitem{serre1}  D. Serre,  Solutions faibles globales des quations de Navier-Stokes pour un fluide compressible. \emph{C. R. Acad. Sci. Paris Sér. I Math.} 303(13)(1986), 639-642.

    \bibitem{serre2} D. Serre,  On the one-dimensional equation of a viscous, compressible, heat-conducting fluid. \emph{C. R. Acad. Sci. Paris Sér. I Math.} 303(14)(1986), 703-706.

\bibitem{is2002} I. Stra\v{s}kraba, A. Zlotnik, On a decay rate for 1D-viscous compressible barotropic fluid equations, \textit{J. Evolution Equations},    \textbf{2}(2002), 69-96.

 \bibitem{is2003} I. Stra\v{s}kraba, A. Zlotnik, Global properties of solutions to 1D-viscous compressible
barotropic fluid equations with density dependent viscosity, \textit{Z. Angew. Math. Phys.},    \textbf{54}(2003), 593-607.


\bibitem{yyz}T. Yang, Z. A. Yao, C. J. Zhu, Compressible Navier-Stokes equations with density-dependent viscosity and vacuum, \emph{Comm.
Partial Differential Equations}, \textbf{26}(2001), 965-981.

\bibitem{yz}T. Yang, C. J. Zhu, Compressible Navier-Stokes equations with degenerate viscosity coefficient and vacuum, \emph{Comm. Math. Phys.}, \textbf{230}(2002), 329-363.


 \bibitem{zl1} A. A.  Zlotnik,    Uniform estimates and stabilization of symmetric
 solutions of a system of quasilinear equations.  \textit{Diff. Eqs.},  \textbf{36}(2000),  701-716.

\bibitem{zaa95} A. A. Zlotnik,   Uniform estimates and stabilization of solutions to equations of one-dimensional motion
of a multicomponent barotropic mixture. \textit{Math. Notes},
  \textbf{58}(1995),  885-889.

\end {thebibliography}
\end{document}